\documentstyle[12pt,twoside]{article}
\newtheorem{theorem}{\bf Theorem}[section]
\newtheorem{proposition}[theorem]{\bf Proposition}

\newtheorem{corollary}[theorem]{\bf Corollary}

\input{amssym}

\date{}
\begin{document}

\title{{\Large\bf Left $\theta$-derivations on  weighted convolution algebras}}

\author{{\normalsize\sc M. Eisaei, M. J. Mehdipour and Gh. R. Moghimi}}
\maketitle

{\footnotesize  {\bf Abstract.} Let $\theta$ be a homomorphism on $L_0^\infty({\Bbb R}^+, \omega)^*$. In this paper, we study left $\theta$-derivations on $L_0^\infty({\Bbb R}^+, \omega)^*$. We show that every left $\theta$-derivation on $L_0^\infty({\Bbb R}^+, \omega)^*$ is always a $\theta$-derivation, and if $\theta$ is isomorphism, then $L_0^\infty({\Bbb R}^+, \omega)^*$ has no non-zero  left $\theta$-derivation. We also investigate automatic continuity, Singer-Wermer's conjecture and Posner's first theorem for left $\theta-$derivations on $L_0^\infty({\Bbb R}^+, \omega)^*$.}
{
}
\section{\normalsize\bf Introduction}

Throughout this paper, $\omega$ is a weight function on ${\Bbb R}^ += [ { 0 ,\infty })$. Let us recall that a \emph{weight function} on ${\Bbb R}^ + $ is a continuous function
$ \omega: {\Bbb R}^ +\rightarrow[ {1, + \infty }) $ such that $ \omega \left(  0  \right) = 1 $ and for every $ x,y \in{\Bbb R}^ + $
\[\omega \left( {x + y} \right) \leqslant \omega \left( x \right)\omega \left( y \right).\]
Let $  L^1({\Bbb R}^+, \omega)  $ be the Banach space of all Lebesgue measurable functions $ f $ on ${\Bbb R}^ +$ with $ \omega f \in {L^1}({\Bbb R}^ +) $, the Banach algebra of all Lebesgue integrable functions  on ${\Bbb R}^ +$. 
Let also $ M({\Bbb R}^+, \omega) $ be the Banach algebra of all complex regular Borel measures $ \mu $ on $ {\Bbb R}^ + $ such that $ \omega \mu  \in M( {\Bbb R}^ +) $, the Banach algebra of all complex regular Borel measures on $ {\Bbb R}^ + $. It is well-known that $M({\Bbb R}^+, \omega)  $ 
is a commutative Banach algebra with the convolution product and the norm 
$$
\|\mu\|_\omega=\|\omega\mu\|\quad\quad(\mu\in M({\Bbb R}^+, \omega)).
$$
Also, $ L^1({\Bbb R}^+, \omega)  $ can be regarded as a closed ideal in $ M({\Bbb R}^+, \omega)$; for a study of these spaces see [7, 14].

 Let $ {L_0^\infty({\Bbb R}^+, \omega)} $ be the Banach space of all Lebesgue measurable functions $ f $ on $ {\Bbb R}^ + $ such that
$$
\lim_{x\rightarrow\infty}\hbox{ess sup}\{\frac{f(y)\chi_{(x, \infty)}(y)}{\omega(y)}: y\in{\Bbb R}^+\}=0.
$$
It is proved that the dual of $ L_0^\infty({\Bbb R}^+, \omega)$ is a Banach algebra with the first Arens product
defined by
$
\langle {s\cdot t, f}\rangle= \langle {s, t f} \rangle,
$
where $\langle {t f, g}\rangle=\langle {t, fg} \rangle$ and 
$$
fg(x)= \int_ 0 ^\infty  {f\left( {x + y} \right)g\left( y \right)dy}
$$
for all $ s, t \in  L_0^\infty({\Bbb R}^+, \omega)^*$, $f \in L_0^\infty({\Bbb R}^+, \omega)^* $ and $g \in  L^1({\Bbb R}^+, \omega)  $; for more details see [9, 10, 11]. We denote by $ \Lambda \left( L_0^\infty({\Bbb R}^+, \omega)^* \right) $ the set of all right identities of $ L_0^\infty({\Bbb R}^+, \omega)^* $ with norm one. Let us remark that if $ \nu\in \Lambda \left( L_0^\infty({\Bbb R}^+, \omega) ^*\right) $, then $ \nu\cdot L_0^\infty({\Bbb R}^+, \omega)^* $ and $ M \left({\Bbb R}^+, \omega  \right) $ are isometric isomorphic. Since $ M({\Bbb R}^+, \omega)  $ is a commutative Banach algebra, $ \nu \cdot L_0^\infty({\Bbb R}^+, \omega)^* $ is commutative and hence
$$
k\cdot   s\cdot t = k\cdot\nu \cdot   s\cdot    \nu \cdot   t = k\cdot   \nu \cdot   t\cdot   \nu \cdot   s = k\cdot   t\cdot   s
$$
for all $ k, s, t\in L_0^\infty({\Bbb R}^+, \omega)^* $. By Titchmarsh's convolution theorem [3], $ M({\Bbb R}^+, \omega)  $ is an integral domain. So $ \nu \cdot L_0^\infty({\Bbb R}^+, \omega)^* $ is an integral domain; see also Proposition 5.1 of [1]. We denote the right annihilator $L_0^\infty({\Bbb R}^+, \omega)^*$ by $\hbox{ran}(L_0^\infty({\Bbb R}^+, \omega)^*) $ and defind by
$$
\hbox{ran}(L_0^\infty({\Bbb R}^+, \omega)^*)=\{r \in L_0^\infty({\Bbb R}^+, \omega)^*: s\cdot    r=  0\quad\hbox{for all}\quad s\in  L_0^\infty({\Bbb R}^+, \omega)^*\}.
$$
It is easy to see that 
$$
L_0^\infty({\Bbb R}^+, \omega)^* = \nu \cdot   L_0^\infty({\Bbb R}^+, \omega)^* \oplus\hbox{ran}(L_0^\infty({\Bbb R}^+, \omega)^*).
$$
Let $\theta$ be a homomorphism on $L_0^\infty({\Bbb R}^+, \omega)^*$. A linear operator $\delta: L_0^\infty({\Bbb R}^+, \omega)^*\rightarrow L_0^\infty({\Bbb R}^+, \omega)^*$ is called a \emph{left $\theta-$derivation} if
$$
\delta(s\cdot t)=\theta(s)\cdot \delta(t)+\theta(t)\cdot \delta(s)
$$
for all $s, t\in  L_0^\infty({\Bbb R}^+, \omega)^*$. Also, $\delta$ is a \emph{$\theta-$derivation} if
$$
\delta(s\cdot t)=\theta(s)\cdot \delta(t)+\delta(s)\cdot \theta(t)
$$
for all $s, t\in  L_0^\infty({\Bbb R}^+, \omega)^*$.

Ashraf [3] studied Jordan left $\theta$-derivations and left $\theta$-derivations of a prime ring $R$ with characteristic different from 2. He proved that if $d$ is a Jordan left $\theta$-derivation of $R$, then $d=0$ or $R$ is commutative; see [5, 8] for Jordan left derivations. He also showed that every Jordan left $\theta$-derivation is a left $\theta$-derivation; see [4] for generalized left $\theta$-derivations.   Here, we investigate the results for Banach algebra $ L_0^\infty({\Bbb R}^+, \omega)^*$. Note that it is easy to see that 
$ L_0^\infty({\Bbb R}^+, \omega)^*$ is not a prime ring.

In this paper, we investigate left $ \theta $-derivations on $ L_0^\infty({\Bbb R}^+, \omega)^*$ and show that every left $\theta$-derivation on  $ L_0^\infty({\Bbb R}^+, \omega)^*$ is a  $\theta$-derivation on  $ L_0^\infty({\Bbb R}^+, \omega)^*$. In the case where $\theta$ is isomorphism, the zero map is the only left $\theta$-derivation on  $ L_0^\infty({\Bbb R}^+, \omega)^*$. We also study continuity of left $\theta$-derivations and prove that if $\|\theta(\nu )\|<1$, then the only continuous left $\theta$-derivation on $L_0^\infty({\Bbb R}^+, \omega)^*$ is zero, where $\nu \in\Lambda(L_0^\infty({\Bbb R}^+, \omega)^*)$. Then we give an analog of Posner's first theorem. Finally, we investigate Singer-Wermer's conjecture and establish that a left $\theta-$derivation on $ L_0^\infty({\Bbb R}^+, \omega)^*$ maps $ L_0^\infty({\Bbb R}^+, \omega)^*$ into its radical if and only if it is zero. 

\section{\normalsize\bf Left $ \theta $-derivations on $ L_0^\infty({\Bbb R}^+, \omega)^*$}

First, we investigate homomorphisms on $ L_0^\infty({\Bbb R}^+, \omega)^*$.

\begin{proposition}\label{0-1} 
Let $ \theta $ be a non-zero homomorphism on $ L_0^\infty({\Bbb R}^+, \omega)^*$. Then the following statements hold.

\emph{(i)} If $ \nu \in \Lambda(L_0^\infty({\Bbb R}^+, \omega)^*)$, then 
$ \theta(\nu ) = \nu + r_{0} $ for some $ {r_ 0 } \in\emph{ran}(L_0^\infty({\Bbb R}^+, \omega)^*) $.

\emph{(ii)} $ \theta ( \emph{ran}(L_0^\infty({\Bbb R}^+, \omega)^*) ) \subseteq \emph{ran}(L_0^\infty({\Bbb R}^+, \omega)^*) ) $.

\emph{(iii)} If  $ \theta $ is monomorphism, then $ \theta ( \emph{ran}(L_0^\infty({\Bbb R}^+, \omega)^*) )=\emph{ran}(L_0^\infty({\Bbb R}^+, \omega)^*)\cap\emph{Im}(\theta)$.

\emph{(iv)} If $ \theta $ is isometric, then 
  $ \theta(\nu )=\nu $ for all $\nu \in \Lambda ( L_0^\infty({\Bbb R}^+, \omega)^*)$. 
\end{proposition}
{\it Proof.}
(i) Let $ \nu \in \Lambda(L_0^\infty({\Bbb R}^+, \omega)^*)$. Then there exist $ {s_0  } \in L_0^\infty({\Bbb R}^+, \omega)^* $ and $r_ 0\in \hbox{ran}( L_0^\infty({\Bbb R}^+, \omega)^*) $ such that $ \theta(\nu ) = \nu \cdot {s_0  } + {r_ 0 } $.  So we have
\begin{eqnarray*}
( {\nu \cdot {s_0  } + {r_ 0 }} )\cdot( {\nu \cdot {s_0  } + {r_ 0 }} )=\theta(\nu )\cdot \theta(\nu )=\theta(\nu ) = \nu \cdot {s_0  } + {r_ 0 }.
\end{eqnarray*}
It follows that
\[( {\nu \cdot {s_0  }} )\cdot( {\nu \cdot {s_0  }} ) - \nu \cdot {s_0  } = {r_ 0 } - {r_ 0 }\cdot {s_0  }.\]
From this and the fact that $$ \nu \cdot L_0^\infty({\Bbb R}^+, \omega)^* \cap\hbox{ran}( L_0^\infty({\Bbb R}^+, \omega)^*) = \{  0  \} ,$$ we conclude that
\[( {\nu \cdot {s_0  }} )\cdot( {\nu \cdot {s_0  }} ) - \nu \cdot {s_0  } =  0 . \]
Since $ \nu \cdot L_0^\infty({\Bbb R}^+, \omega)^* $ is an integral domain, $ \nu \cdot {s_0  } =  0 $ or $ \nu \cdot {s_0  } =\nu $.
If $ \nu \cdot {s_0  } =  0 $, then $ \theta(\nu ) = {r_ 0 } $ and thus
\[\theta(\nu ) = \theta(\nu )\cdot \theta(\nu ) = {r_ 0 }\cdot {r_ 0 } =  0 .\]
Hence for any $ s\in  L_0^\infty({\Bbb R}^+, \omega)^* $, we have
\[\theta (s) = \theta ( {s\cdot  \nu} ) = \theta (s)\cdot \theta(\nu ) =  0 .\]
That is, $ \theta = 0 $, a contradiction. So $ \nu \cdot s_0 =\nu $ and therefore, 
 $ \theta(\nu ) = \nu + {r_ 0 } $.

(ii) Let $ r \in\hbox{ran}( L_0^\infty({\Bbb R}^+, \omega)^*) $ and $ \nu \in \Lambda(L_0^\infty({\Bbb R}^+, \omega)^*)$. Then
\[ 0  = \theta ( {\nu \cdot r} ) = \theta(\nu )\cdot \theta ( r ) = ( {\nu + {r_ 0 }} )\cdot \theta ( r )\]
for some $ {r_ 0 } \in\hbox{ran}( L_0^\infty({\Bbb R}^+, \omega)^*) $. This shows that
\[s\cdot  \theta ( r ) = s\cdot  ( {\nu + {r_ 0 }} )\cdot \theta ( r ) =  0 \]
for all $ s\in  L_0^\infty({\Bbb R}^+, \omega)^* $. That is, $ \theta ( r ) \in \hbox{ran}( L_0^\infty({\Bbb R}^+, \omega)^*) $.

(iii)  Let $ {r} \in\hbox{ran}( L_0^\infty({\Bbb R}^+, \omega)^*)\cap\hbox{Im}(\theta)$. Then $ \theta(s) = {r} $ for some 
  $s \in L_0^\infty({\Bbb R}^+, \omega)^*$. For every $ t\in  L_0^\infty({\Bbb R}^+, \omega)^* $, we have
 \[\theta ( t\cdot s  ) = \theta(t)\cdot \theta(s) = \theta(t)\cdot {r} =  0 .\]
 Hence $ \theta ( t\cdot s  ) =  0 $. By injectivity of $ \theta $,
 \[t\cdot s  =  0\]
for all $ t\in  L_0^\infty({\Bbb R}^+, \omega)^* $. That is, $  s  \in \hbox{ran}( L_0^\infty({\Bbb R}^+, \omega)^*) $. Therefore,
  $\hbox{ran}( L_0^\infty({\Bbb R}^+, \omega)^*)\cap\hbox{Im}(\theta)$ is a subset of $\theta ( {\hbox{ran}( L_0^\infty({\Bbb R}^+, \omega)^*)})$.

(iv)  Let $ \nu \in \Lambda(L_0^\infty({\Bbb R}^+, \omega)^*)$. Then $ \theta(\nu ) = \nu + {r_ 0 } $ for some 
  $ {r_ 0 } \in \hbox{ran}( L_0^\infty({\Bbb R}^+, \omega)^*) $. For every $ s\in  L_0^\infty({\Bbb R}^+, \omega)^* $ we have $ s\cdot  \theta(\nu ) = s $. Thus $ \theta(\nu ) $ is a right identity for $ L_0^\infty({\Bbb R}^+, \omega)^*$. Since $ \theta $ is isometry, 
  $ \| {\theta(\nu )} \| = \|\nu \| = 1 $. So $ \theta(\nu )=\nu _1 \in \Lambda(L_0^\infty({\Bbb R}^+, \omega)^*)$ and consequently, $\nu -\nu _1=r_0$. Hence, $\nu =\nu _1$. 
 $\hfill\square$\\

For a homomorphism $\theta$ on $L_0^\infty({\Bbb R}^+, \omega)^*$, let us recall that a left $ \theta $-derivations $\delta$ on $ L_0^\infty({\Bbb R}^+, \omega)^*$ is called $\theta$-\emph{commuting} if for every $s\in  L_0^\infty({\Bbb R}^+, \omega)^*$
$$
\delta(s)\cdot\theta(s)=\theta(s)\cdot\delta(s).
$$
In the case where $\theta=id$, $\theta$
is called \emph{commuting}, where $id$ is the identity map on $L_0^\infty({\Bbb R}^+, \omega)^*$. Now, we are ready to prove the main result of this paper.

\begin{theorem}\label{0-4} 
Let  $ \theta $ be a homomorphism on $ L_0^\infty({\Bbb R}^+, \omega)^*$ and $ \nu \in \Lambda(L_0^\infty({\Bbb R}^+, \omega)^*)$. If $ \delta $ is a left $ \theta $-derivation on $ L_0^\infty({\Bbb R}^+, \omega)^*$, then the following statements hold.

 \emph{(i)} $ \delta( \emph{ran}(L_0^\infty({\Bbb R}^+, \omega)^*) ) = \delta ( {\Lambda ( L_0^\infty({\Bbb R}^+, \omega) ^*)} ) = \{  0  \}$.

 \emph{(ii)} $\delta(s)=\theta(\nu )\cdot\delta(s)$ for all $s\in  L_0^\infty({\Bbb R}^+, \omega)^* $.

 \emph{(iii)} $\delta(s)\cdot\theta(t)=\theta(t)\cdot\delta(s)$ for all $s, t\in  L_0^\infty({\Bbb R}^+, \omega)^*$. 

\emph{(iv)} $\delta $ is a $ \theta $-commating $ \theta $-derivation.
\end{theorem}
{\it Proof.} (i) Let $ \nu \in \Lambda(L_0^\infty({\Bbb R}^+, \omega)^*)$. Then
\begin{equation}\label{1} 
\delta(\nu ) = \delta ( {\nu \cdot\nu} ) = \theta(\nu )\cdot \delta(\nu ) + \theta(\nu )\cdot \delta(\nu ) = 2\theta(\nu )\cdot \delta(\nu ).
\end{equation}
It follows that
\[\delta(\nu ) = 4\theta(\nu )\cdot \delta(\nu ).\]
Therefore, $ \theta(\nu )\cdot \delta(\nu ) =  0 $. According to (1), $ \delta(\nu ) =  0$. 

Assume now that $ r \in\hbox{ran}( L_0^\infty({\Bbb R}^+, \omega)^*) $. Then by (1) we have
\begin{eqnarray}\label{par}
  0  &=& \delta( \nu\cdot r)\nonumber\\
&=& \theta (\nu)\cdot\delta ( r ) + \theta ( r )\cdot\delta (\nu)\nonumber\\
&=& \theta (\nu)\cdot\delta ( r ) + \theta ( r )\cdot( {2\theta (\nu)\cdot\delta ( \nu )} ) \\
&=& \theta (\nu)\cdot\delta ( r ) + 2\theta ( r )\cdot\delta (\nu )\nonumber\\
&=&\nu\cdot\delta(r)+r_0\cdot\delta(r)+2\theta ( r )\cdot\delta (\nu )\nonumber,
\end{eqnarray}
where $\theta(\nu)=\nu+r_0$ for some $r_0\in \hbox{ran}( L_0^\infty({\Bbb R}^+, \omega)^*) $. This together with Proposition 2.1(ii) and the fact that 
$$
\nu\cdot L_0^\infty({\Bbb R}^+, \omega)^* \cap\hbox{ran}( L_0^\infty({\Bbb R}^+, \omega)^*)=\{0\}$$
shows that $\nu\cdot\delta(r)=\{0\}$ and hence $\delta(r)\in \hbox{ran}( L_0^\infty({\Bbb R}^+, \omega)^*)$. Thus  $\theta (\nu)\cdot\delta ( r ) =0$. From this and (2) we see that 
 $ \theta ( r )\cdot\delta (\nu )=  0 $. So
\[\delta ( r ) = \delta ( {r\cdot\nu} ) = \theta ( r )\cdot \delta(\nu ) + \theta(\nu )\cdot \delta ( r ) =  0 .\]
That is, $ \delta$ is zero on ${\hbox{ran}( {L_ 0 ^\infty {{( {\Bbb R}^+, \omega )}^ * }} )}$. 

(ii) Let $s\in L_0^\infty({\Bbb R}^+, \omega)^*$ and $\nu \in \Lambda(L_0^\infty({\Bbb R}^+, \omega)^*)$. Then $s-\nu\cdot s\in\hbox{ran}(L_0^\infty({\Bbb R}^+, \omega)^*)$. By (i), we have
$\delta(s-\nu\cdot s)=0$. Since $\delta(\nu)=0$, we have
$$
\delta(s)=\delta(\nu\cdot s)=\theta(\nu)\cdot \delta(s)+ \theta(s)\cdot\delta(\nu)=\theta(\nu)\cdot \delta(s).
$$
(iii) Let $ \delta $ be a left $ \theta $-derivation, $\nu \in \Lambda(L_0^\infty({\Bbb R}^+, \omega)^*)$, $ s, t \in  L_0^\infty({\Bbb R}^+, \omega)^* $ and $ r \in \hbox{ran}( L_0^\infty({\Bbb R}^+, \omega)^*) $. Then
\[ 0  = \delta ( {s\cdot  r} ) = \theta (s)\cdot \delta ( r ) + \theta ( r )\cdot \delta(s) = \theta ( r )\cdot \delta(s)\]
and
\[\delta(s)\cdot \theta ( {\nu \cdot t} ) = \theta(\nu )\cdot \delta(s)\cdot \theta ( {\nu \cdot t} ) = \theta ( {\nu \cdot t} )\cdot \delta(s).\]
Choose $ {r_s},{r_t} \in\hbox{ran}( L_0^\infty({\Bbb R}^+, \omega)^*) $ such that
\[s = \nu \cdot s + {r_s}\quad\hbox{and}\quad t = \nu \cdot t + {r_t}.\]
Since $ \delta ( {\hbox{ran}( L_0^\infty({\Bbb R}^+, \omega)^*)} ) = \{  0  \} $, it follows that $\theta ( {\nu \cdot s} )\cdot \delta ( {{r_t}} )=  0  $ and so
\[\theta ( {{r_t}} )\cdot \delta ( {\nu \cdot s} ) = \delta ( {( {{\nu \cdot s}} )\cdot {r_t}} ) =  0 .\]
We also have
\begin{eqnarray*}
\delta ( {\nu \cdot s} )\cdot \theta ( {\nu \cdot t} )&=& \theta(\nu )\cdot \delta(s)\cdot \theta ( {\nu \cdot t} )\\
&=& \theta ( {\nu \cdot t} )\cdot \delta(s)\\
&=& \theta ( {\nu \cdot t} )\cdot \delta ( {\nu \cdot s} ).
\end{eqnarray*}
These imply that $ \theta(t)\cdot \delta(s) = \delta(s)\cdot \theta(t) $. 

(iv) This follows at once from (iii).
$\hfill\square$\\

We now give some consequence of Theorem 2.2.

\begin{corollary}
Let  $ \theta $ be a  homomorphism on $ L_0^\infty({\Bbb R}^+, \omega)^*$. If $\delta $ is a left $ \theta $-derivation on $ L_0^\infty({\Bbb R}^+, \omega)^*$, then $\delta = 0 $.
\end{corollary}
{\it Proof.}
Let  $ \nu \in \Lambda(L_0^\infty({\Bbb R}^+, \omega)^*)$. Then $\delta(\nu)=0$ and so
\[\delta(s) = \delta( {s\cdot\nu} ) = \delta(s)\cdot \delta(\nu ) =  0\]
for all $ s\in  L_0^\infty({\Bbb R}^+, \omega)^*$.
 $\hfill\square$

\begin{corollary}\label{zmz} Let  $ \theta $ be a  homomorphism on $ L_0^\infty({\Bbb R}^+, \omega)^*$, $\delta$ be a non-zero left $ \theta $-derivation on $ L_0^\infty({\Bbb R}^+, \omega)^*$ and $\nu \in\Lambda(L_0^\infty({\Bbb R}^+, \omega)^*)$. Then $\delta$ is continuous if and only if $\delta|_{\nu \cdot L_0^\infty({\Bbb R}^+, \omega)^* }$ is continuous. In this case, $\|\theta(\nu )\|\geq 1$.
\end{corollary}
{\it Proof.} Let $\delta|_{\nu\cdot L_0^\infty({\Bbb R}^+, \omega)^* }$ be continuous.
If $\{s_i\}_i$ is a sequence in $L_0^\infty({\Bbb R}^+, \omega)^*$ with 
$$
s_i\rightarrow s\in  L_0^\infty({\Bbb R}^+, \omega)^*,
$$
then $\nu \cdot s_i\rightarrow \nu \cdot s.$ It follows that
$$
\delta(s_i)=\delta(\nu \cdot s_i)\rightarrow\delta(\nu \cdot s)=\delta(s).
$$
Hence $\delta$ is continuous. The converse is clear.

To complete the proof, let $s\in  L_0^\infty({\Bbb R}^+, \omega)^*$ such that $\delta(s)\neq 0$. Then
$$
\|\delta(s)\|=\|\theta(\nu )\cdot\delta(s)\|\leq\|\theta(\nu )\|\|\delta(s)\|.
$$
Hence $\|\theta(\nu )\|\geq 1$. $\hfill\square$\\

As an immediate consequence of Corollary 2.4,  we have the following result.

\begin{corollary} Let $\nu \in\Lambda(L_0^\infty({\Bbb R}^+, \omega)^*)$ and  $ \theta $ be a  homomorphism on $ L_0^\infty({\Bbb R}^+, \omega)^*$ with $\|\theta(\nu )\|<1$. Then the zero map is the only continuous left $\theta$-derivation on $ L_0^\infty({\Bbb R}^+, \omega)^*$.
\end{corollary}

Note that if $\theta$ is a monomorphism on $ L_0^\infty({\Bbb R}^+, \omega)^*$, then there exists a map $\bar{\theta}$ on $ L_0^\infty({\Bbb R}^+, \omega)^*$ such that $\bar{\theta}\theta=id$.

\begin{theorem} 
Let  $ \theta $ be a homomorphism on $ L_0^\infty({\Bbb R}^+, \omega)^*$. If $ \delta $ is a left $ \theta $-derivation on $ L_0^\infty({\Bbb R}^+, \omega)^*$, then the following statements hold.

\emph{(i)} $\emph{Im}(\delta)\subseteq\cap\{u\cdot L_0^\infty({\Bbb R}^+, \omega)^*: u\; \emph{is a right identity of}\;   L_0^\infty({\Bbb R}^+, \omega)^*\}$.

\emph{(ii)} If $ \theta $ is monomorphism and $\bar{\theta}$ is homomorphism, then $\emph{Im}(\delta)\subseteq\emph{ker}(\bar{\theta})$ and $\hbox{Im}(\delta)\cap\hbox{Im}(\theta)=\{0\}$. In this case, $\delta$ is not surjective.

\emph{(iii)} If $ \theta $ is isomorphism, then $\delta=0$ on $L_0^\infty({\Bbb R}^+, \omega)^*$.
\end{theorem}
{\it Proof.} (i) Let $\nu\in\Lambda( L_0^\infty({\Bbb R}^+, \omega)^*)$. Then $\theta(\nu)$ is a right identity of $ L_0^\infty({\Bbb R}^+, \omega)^*$ by Proposition 2.1(i). Now, apply Theorem 2.2(ii).

(ii) Let  $ \theta $ be a monomorphism on $ L_0^\infty({\Bbb R}^+, \omega)^*$.  Set $\eta:=\bar{\theta}\delta$. Then $\eta$ is a left derivation on $ L_0^\infty({\Bbb R}^+, \omega)^*$ and so by Theorem 4.1 of [1], $$\bar{\theta}\delta=\eta=0.$$ This implies that $$\hbox{Im}(\delta)\subseteq\hbox{ker}(\bar{\theta}).$$
Let $s, t\in  L_0^\infty({\Bbb R}^+, \omega)^*$ and $\delta(s)=\theta(t)$. Then
$$
0=\bar{\theta}(\delta)(s)=\bar{\theta}\theta(t)=t.
$$
This shows that $\hbox{Im}(\delta)\cap\hbox{Im}(\theta)=\{0\}$. That is, (i) holds.  Finally, suppose that $\delta$ is surjective. Then
$$
L_0^\infty({\Bbb R}^+, \omega)^*=\hbox{Im}(\delta)\subseteq\hbox{ker}(\bar{\theta}).
$$
Therefore, $\bar{\theta}=0$, a contradiction.

(iii) Let $\theta$ be isomorphism. Then $\hbox{ker}(\bar{\theta})=0$. So by (ii), $\delta$ is zero on $L_0^\infty({\Bbb R}^+, \omega)^*$.
$\hfill\square$\\

Singer-Wermer's conjecture states that the range of a derivation of a Banach algebra $A$ is into $\hbox{rad}(A)$, the radical of $A$. Bre\v{s}ar and Vukmann [5] studied the conjecture for left derivations and proved that the range of a continuous left derivation of a Banach algebra is into its radical. Now, we investigate the conjecture for left $\theta$-derivations of $ L_0^\infty({\Bbb R}^+, \omega)^*$. 

\begin{theorem}\label{sw} Let  $ \theta $ be a monomorphism on $ L_0^\infty({\Bbb R}^+, \omega)^*$ and $\delta$ be a left $\theta$-derivation on $ L_0^\infty({\Bbb R}^+, \omega)^*$.  Then the range of $\delta$ is contained into $\emph{rad}(L_0^\infty({\Bbb R}^+, \omega)^*)$ if and only if $\delta=0$.
\end{theorem}
{\it Proof.} Let $s\in L_0^\infty({\Bbb R}^+, \omega)^*$ and $\delta(s)\in\hbox{rad}(L_0^\infty({\Bbb R}^+, \omega)^*)$. By Theorem 2.1 of [12], we have
$$
\hbox{rad}(L_0^\infty({\Bbb R}^+, \omega)^*)=\hbox{ran}(L_0^\infty({\Bbb R}^+, \omega)^*).
$$
Hence $\delta(s)\in\hbox{ran}(L_0^\infty({\Bbb R}^+, \omega)^*)$. In view of Theorem 2.2(ii), we have
$$
\delta(s)=\theta(\nu)\cdot\delta(s)=0,
$$
where $\nu\in\Lambda(L_0^\infty({\Bbb R}^+, \omega)^*)$.
Therefore, $\delta=0.$
$\hfill\square$\\

Posenr [13] proved that if the product of two derivations of a prime ring of characteristic not 2 is a derivation, then the one of derivations is zero; see Creedon [6] for product of derivations of Banach algebras. This result is known as Posner's first theorem.

\begin{theorem} Let  $ \theta $ be an idempotent monomorphism on $ L_0^\infty({\Bbb R}^+, \omega)^*$, and let $\delta_1$, $\delta_2$ and $\delta_1\delta_2$ be left $\theta$-derivations on $ L_0^\infty({\Bbb R}^+, \omega)^*$. Then the following statements hold.

\emph{(i)} $\delta_1|_{\emph{Im}(\theta)}=0$ or $\delta_2=0$.

\emph{(ii)} If $\theta\delta_1=\delta_1\theta$, then $\delta_1=0$ or $\delta_2=0$.
\end{theorem}
{\it Proof.} Since $\delta_1$, $\delta_2$ and $\delta_1\delta_2$ are left $\theta$-derivations on $ L_0^\infty({\Bbb R}^+, \omega)^*$, it follows that 
$$
\theta\delta_2(s)\delta_1\theta(s)=0
$$
for all $s\in L_0^\infty({\Bbb R}^+, \omega)^*$. In view of Proposition 5.1 of [1], we have
$$
\theta\delta_2(s)\in\hbox{ran}(L_0^\infty({\Bbb R}^+, \omega)^*)\quad\hbox{or}\quad\delta_1\theta(s)\in \hbox{ran}(L_0^\infty({\Bbb R}^+, \omega)^*).
$$
If $\theta\delta_2(s)\in\hbox{ran}(L_0^\infty({\Bbb R}^+, \omega)^*)$, then Proposition 2.1 (iii) implies that $\delta_2(s)\in\hbox{ran}(L_0^\infty({\Bbb R}^+, \omega)^*)$. By Theorem 2.7, $\delta_2=0$. If $\delta_1\theta(s)\in \hbox{ran}(L_0^\infty({\Bbb R}^+, \omega)^*)$, then another application of Theorem 2.7 shows that $\delta_1\theta(s)=0$ for all $s\in L_0^\infty({\Bbb R}^+, \omega)^*$. Therefore, $\delta_1|_{\hbox{Im}(\theta)}=0$. That is, (i) holds. For (ii), we only recall that if $\theta\delta_1=\delta_1\theta$, then by the proof of (i)
$$
\theta\delta_1=\delta_1\theta=0.
$$
This implies that $\delta_1=0$ on $L_0^\infty({\Bbb R}^+, \omega)^*$.
$\hfill\square$

\begin{corollary} Let  $ \theta $ be an idempotent monomorphism on $ L_0^\infty({\Bbb R}^+, \omega)^*$, and let $\delta$ and $\delta^2$ be left $\theta$-derivations on $ L_0^\infty({\Bbb R}^+, \omega)^*$. If $\delta\theta=\theta\delta$, then $\delta$ is zero.
\end{corollary}


\vspace{2mm}
 {\footnotesize

\noindent {\bf Mojdeh Eisaei}\\
Department of Mathematics,\\ Payame Noor University (PNU),\\
Tehran, Iran\\ e-mail: mojdehessaei59@student.pnu.ac.ir\\
{\bf Mohammad Javad Mehdipour}\\
Department of Mathematics,\\ Shiraz University of Technology,\\
Shiraz
71555-313, Iran\\ e-mail: mehdipour@sutech.ac.ir\\
{\bf Gholam Reza Moghimi}\\
Department of Mathematics,\\
Payame Noor University (PNU),\\
Tehran, Iran\\
e-mail: moghimimath@pnu.ac.ir
\end{document}